\documentclass[10pt]{amsart}

\usepackage{amsmath}
\usepackage{amsthm}
\usepackage{amsfonts}
\usepackage{amssymb,latexsym}
\usepackage[all]{xy}
\usepackage{graphicx}
\usepackage[mathscr]{eucal}
\usepackage{verbatim}

\theoremstyle{plain}
    \newtheorem{thm}{Theorem}[section]
    \newtheorem{prop}[thm]{Proposition}
    \newtheorem{lemma}[thm]{Lemma}
    
    \newtheorem{cor}[thm]{Corollary}

\theoremstyle{definition}
    \newtheorem{defn}[thm]{Definition}

\theoremstyle{remark}
    \newtheorem{rem}[thm]{Remark}
    \newtheorem{example}[thm]{Example}

\newtheorem*{namedtheorem}{\theoremname}
\newcommand{\theoremname}{testing}

\newcommand{\A}{\mathcal A}

\newcommand{\C}{\mathcal C}

\newcommand{\T}{\mathcal T}
\newcommand{\W}{\mathcal W}

\newcommand{\stk}[1]{\stackrel{#1}{\longrightarrow}}

\newcommand{\ints}{\ensuremath{\mathbb{Z}}}

\newcommand{\rar}{\ensuremath{\rightarrow}}

\DeclareMathOperator{\Supp}{Supp}
\DeclareMathOperator{\Spec}{Spec}
\DeclareMathOperator{\Hom}{Hom}

\DeclareMathOperator{\proj}{proj}

\DeclareMathOperator{\Wide}{Wide}
\DeclareMathOperator{\Thick}{Thick}

\DeclareMathOperator{\Ext}{Ext}

\DeclareMathOperator{\plaincolim}{colim}

\begin{document} 

\title[Abelian subcategories closed under extensions]{Abelian subcategories closed under extensions: \\ $K$-theory and decompositions}
\date{\today}

\author{Sunil K. Chebolu} 
\address {Department of Mathematics \\
          University of Western Ontario \\
          London, ON,  N6A 5B7.}
\email{schebolu@uwo.ca}


\begin{abstract} A full subcategory of modules over a commutative ring $R$ 
is  \emph{wide} if it is  abelian and closed under extensions. Hovey ~\cite{wide} gave a 
classification of the wide subcategories of finitely presented modules over regular coherent rings in terms of
certain specialisation closed subsets of $\Spec(R)$. We use this classification theorem to study $K$-theory and 
Krull-Schmidt type decompositions for wide subcategories. It is shown that the $K$-group, in the sense of Grothendieck, of a wide 
subcategory $\W$ of finitely presented modules over a regular coherent ring is isomorphic to that of the thick subcategory of 
perfect complexes whose homology groups belong to $\W$. We also show that  the wide subcategories of finitely generated modules
over a noetherian regular ring can be decomposed uniquely into indecomposable ones. This result is then applied to obtain 
a decomposition for the $K$-groups of wide subcategories.

\end{abstract}

\maketitle


\section{Introduction}

It is well known that the derived category $D(R)$ of commutative ring $R$ is formally very similar to the stable homotopy category of spectra. For instance, both 
these categories have the structure of a tensor triangulated category. Exploring the deeper structural similarity between these algebraic and topological worlds has led to 
some active research in the last twenty years; see \cite{mps} for a systematic presentation of this research.
This line of research was initiated by Hopkins in the 1980's. In his influential paper \cite{Ho}, Hopkins 
gave an elegant classification of the thick subcategories (triangulated subcategories that are closed under direct summands) of finite spectra and those of perfect complexes 
in $D(R)$. These  classification results have had tremendous impact in their respective fields. Since then there has been a lot of interest in the classification of 
thick subcategories in various other fields including algebraic geometry \cite{Th} and modular representation theory \cite{bcr}. 

The celebrated thick subcategory theorem in $D(R)$ states that the thick subcategories of small objects correspond to subsets of $\Spec(R)$ that are a 
union of subsets $V(I)$, where $I$ is a finitely generated ideal of $R$. Such subsets of $\Spec(R)$ will be called \emph{thick supports}.
In particular, when $R$ is noetherian, the thick subcategories of small objects correspond to arbitrary union of closed subsets of $\Spec(R)$.
(Although this result is due to Hopkins, his original proof has a mistake which was corrected by Neeman \cite{Ne} in the noetherian case, and was generalised
later to arbitrary commutative rings and schemes by Thomason \cite{Th}.)

Hovey \cite{wide} investigated the analogue of such a classification for the category of $R$-modules. He showed that for some large class of commutative rings
the subcategories of $R$-modules that correspond to thick supports of $\Spec(R)$ are those abelian subcategories of finitely presented $R$-modules 
that are closed under extensions  (see theorem \ref{thm:hoveystheorem} for the precise statement). Following Hovey we will call such subcategories as 
wide subcategories; see definition \ref{defn:wide}. Hovey's proof is homological in nature, relating a wide subcategory of modules to a thick subcategory of complexes. 
In particular, his proof relies on the thick subcategory theorem of the derived category stated in the previous paragraph.

In this paper we build on the work of Hovey and study wide subcategories of $R$-modules. We begin, in the next section, by setting up our  apparatus of rings, 
modules and complexes. We then discuss  Hovey's classification of the wide subcategories (theorem \ref{thm:hoveystheorem}). 
In the last two sections we use this classification to study two things about wide subcategories: $K$-theory in section \ref{sec:ktheorywide} and 
Krull-Schmidt type decompositions in section \ref{sec:KSwide}. In more detail,
in section \ref{sec:ktheorywide} we show that when $R$ is a regular coherent ring, the $K$-group of a wide subcategory $\W$ of finitely presented $R$-modules 
agrees with that of the  thick subcategory $\C$ consisting of perfect complexes whose homology  groups belong to $\W$ (theorem \ref{thm:K_0wide}). 
In section \ref{sec:KSwide} we show that when $R$ is a noetherian regular ring,  the wide  subcategories of finitely generated $R$-modules
can be decomposed uniquely into indecomposable wide subcategories (theorem \ref{thm:KSwide}). This  result is then applied to obtain a decomposition of the $K$-groups 
of wide subcategories. As another application, we give a categorical characterisation of local rings among  noetherian regular rings.

All rings in this paper are assumed to be unital and commutative, and all our subcategories will be full. 
We work with chain complexes $C_*$, rather than cochain complexes,  so the differential lowers the degree by one, $d: C_n \rar C_{n-1}$.  

\medskip
\noindent
\textbf{Acknowledgements:} I would like to thank John Palmieri for helping me with the proof of proposition \ref{prop:0}, and James Zhang 
for some interesting comments on decompositions of subcategories. I also want to thank the anonymous referee for giving some valuable suggestions which helped
me improve the presentation of this paper.

\section{Hovey's classification of wide subcategories}

Here and elsewhere, $R$ will denote a unital and commutative ring.

\begin{defn} \label{defn:wide} \cite{wide}  A full subcategory $\C$ of the category of $R$-modules, or any abelian category,
is said to be \emph{wide} if it is a non-empty  abelian subcategory that is closed under extensions. 
\end{defn}

\begin{example} Easy exercises in algebra tell us that the following are wide subcategories.

\begin{enumerate}
\item Rational vector spaces form a wide subcategory of abelian groups.
\item For every prime $p$, $p$-local abelian groups (groups $G$ such that $G \otimes \mathbb{Z}_p \cong G$) form a wide subcategory of abelian groups.
\item Finitely generated modules over a noetherian ring form a wide subcategory of $R$-modules.
\item Finitely presented modules over a coherent ring (ring for which every finitely generated ideal is finitely presented) form a wide subcategory of $R$-modules.
\end{enumerate}
\end{example}

It is natural to ask for a classification of all the wide subcategories of the abelian category of $R$-modules. Although this question is extremely hard to answer in
this generality, it becomes more accessible, as Hovey shows,  if ones restricts to modules with good finiteness conditions and rings with good homological properties -- 
the modules with good finiteness conditions being finitely presented ones, and the rings with good homological properties being ones that are regular and coherent.
We put these objects in the proper context in the next two paragraphs.

To set the categorical stage for our modules and complexes, consider a smaller subcategory of $R$-modules called $\Wide(R)$ -- the wide subcategory generated 
by the ring $R$ in the category of $R$-modules. 
In other words, $\Wide(R)$ is the intersection of all the wide subcategories of $R$-modules that contain $R$.  
Again, it is hard to characterise $\Wide(R)$ for a general ring $R$. However, when 
$R$ is coherent, it turns out that  $\Wide(R)$ is precisely the category of finitely presented $R$-modules; see \cite[Lemma 1.6]{wide}. 
(Similarly, it is an easy exercise to show that if $R$ is noetherian, then $\Wide(R)$ 
is  the category of all finitely generated $R$-modules.) We therefore focus our attention to the wide subcategories of $\Wide(R)$, and assume that $R$ is at least
coherent. Under the partial order given by the inclusion of subcategories, the wide subcategories of $\Wide(R)$ form a lattice which we denote by $\mathcal{L}_{\Wide} (R)$. 
Analogously, in the 
derived category the  small objects are those that belong to  the thick subcategory generated by $R$,  and the full subcategory of these objects 
will be denoted by $\Thick(R)$. It is well-known  that  $\Thick(R)$ is equivalent to the chain homotopy category of perfect 
complexes (bounded complexes of finitely generated projective $R$-modules); see, for instance, \cite[Proposition 9.6]{ch}. The lattice of all the thick 
subcategories of $\Thick(R)$  will be denoted by $\mathcal{L}_{\Thick} (R)$ (as before, the partial order is with respect to inclusion of subcategories).

Now we introduce the class of rings over which we will work. As mentioned in the introduction, Hovey's strategy for classifying the wide subcategories of $\Wide(R)$ is 
to pass to the derived category and use the known classification of the thick subcategories of $\Thick(R)$. In this process he needs a map (\cite[Theorem 3.5]{wide})
\[ \Delta:  \Wide(R) \rar \Thick(R)\]
with the property that for all $M$ in $\Wide(R)$, $H_*(\Delta(M)) \cong M$ as a graded module.
Of course there is an obvious functor taking $\Wide(R)$ into $D(R)$ which has this property -- the one obtained by taking a module and viewing it as a complex in degree zero.
Unless the (finitely presented) module has finite projective dimension, there is no guarantee that this functor will land in $\Thick(R)$ \cite[Proposition 3.4]{wide}.  
This can be ensured if we assume that our coherent ring $R$ is also regular. 
Recall that a commutative ring is \emph{regular} if every finitely generated ideal has finite projective dimension.  
Regular coherent rings have the following desired homological property.
A coherent ring $R$ is regular if and only if every finitely presented $R$-module has finite projective dimension \cite[Theorem 6.2.1]{sarah}. 
So we restrict ourselves to the class of commutative regular coherent rings. These form a large class of commutative rings. Clearly this class includes all 
noetherian regular rings. The following result tells us how to construct more examples of regular coherent rings.

\begin{thm} \cite{sarah} Let $\{ R_{\alpha}\}$ be a diagram of regular coherent rings such that for each arrow $\alpha \rar \beta$, the ring map 
$R_{\alpha} \rar R_{\beta}$ is flat, i.e., $R_{\beta}$ is a flat as a module over $R_{\alpha}$. Then  $\underset{\alpha}{\plaincolim} \, R_{\alpha}$ is a regular 
coherent ring. 
\end{thm}

As a corollary one observes that the polynomial ring in infinitely many variables over a PID is a  non-noetherian example
of a regular coherent ring.

Hovey's classification of wide subcategories is based on some geometric subsets of $\Spec(R)$ which we now define.
A subset of $\Spec(R)$ is a \emph{thick support} if it is a union of (Zariski) closed sets $V(I)$, where $I$ is a finitely generated ideal in $R$. 
In particular when $R$ is noetherian the thick supports are precisely subsets of $\Spec(R)$ that are a union of closed sets (also known as 
\emph{specialisation closed subsets}). The collection of thick supports of $\Spec(R)$ will be denoted by $\mathbb{S}$, and with respect to the partial 
order given by the inclusion of subsets, $\mathbb{S}$ will be viewed as a lattice.
Finally we recall the definition of support for 
modules and complexes which plays a key role in these classification theorems. If $M$ is any $R$-module, then 
\[ \Supp(M) := \{ p \in \Spec(R): M \otimes_R R_p \ne 0 \}.\]
Similarly if $X$ belongs to $\Thick(R)$, then 
\[ \Supp(X) :=  \{ p \in \Spec(R) : X \otimes_R^L R_p \ne 0\}.\]

Having introduced the necessary terminology and preliminary material, we can now state the main result of \cite{wide}:

\begin{thm} \label{thm:hoveystheorem} \cite{wide}  Let $R$ be regular coherent ring, or more generally the quotient of such a ring by a finitely generated
ideal. Then the lattice $\mathcal{L}_{\Wide}(R)$
is isomorphic to the lattice $\mathbb{S}$ of thick supports of $\Spec(R)$. Under this isomorphism, a wide subcategory
$\W$ corresponds to the thick support $\bigcup_{M\in \W} \Supp(M)$.
\end{thm}

In fact, there is a commutative diagram of lattice bijections:
\begin{figure}[!h]
\[
\xymatrix{
\mathcal{L}_{\Wide}(R)\ar[rrrr]^{f}_{\cong}  &  & & &     \mathcal{L}_{\Thick}(R)  \\
& & \mathbb{S} \ar[ull]^{\xi}_{\cong} \ar[urr]_{\zeta}^{\cong} & & 
}
\]
\caption{Bijection between wide subcategories and thick subcategories}
\label{fig:wideandthick}
\end{figure}

where,
\begin{eqnarray*}
 \xi(A)  & =  &\{ M \in \Wide(R): \Supp(M) \subseteq A \},\\ 
\zeta (A) & = & \{ X \in \Thick(R): \Supp(X)\subseteq A\}, \ \ \text{and }\\ 
f(\W)    &  =  & \{ X \in \Thick(R): H_n(X) \in \W \ \ \text{for all} \ \  n \}. \\
\end{eqnarray*}

As indicated in the above diagram, Hovey arrives at this theorem by passing to the derived category to use Thomason's \cite{Th} thick subcategory theorem (isomorphism 
$\zeta$ in the above diagram which holds for all commutative rings).
In fact, it is in the process of this transition (from the module category to the derived category) that the regularity assumption on the ring enters the picture.

\section{$K$-theory for wide subcategories} \label{sec:ktheorywide}

We begin by reminding the reader the definition of a Grothendieck group of a triangulated category.
Let $\T$ denote a  triangulated category that is
essentially small (i.e., it has only a set of isomorphism classes of
objects). The coproduct in $\T$ will be denoted by $\coprod$, and the suspension by $\Sigma$. 
Then the \emph{Grothendieck group} $ K_0(\T) $ is defined to
be the free abelian group on the isomorphism classes of $\T$ modulo
the Euler relations $[B]=[A]+[C]$, whenever $A\rightarrow B
\rightarrow C \rightarrow \Sigma A $ is an exact triangle in $\T$.
($[X]$ denotes the element in the Grothendieck group  that is
represented by the isomorphism class of the object $X$.) This is
clearly an abelian group with $[0]$ as the identity element and
$[\Sigma X]$ as the inverse of $[X]$. The identity
$[A]+[B]=[A \amalg B]$ holds in the Grothendieck group. Also note that any
element of $K_0(\T)$ is of the form $[X]$ for some $X\in \T$. All
these facts follow easily from the axioms for a triangulated category. 
The Grothendieck group of an essentially small abelian category is defined similarly where the relations 
come from short exact sequences in the category.

Since a wide subcategory is also an abelian category, it makes sense to talk
of its Grothendieck group. So we now try to relate the Grothendieck groups of wide subcategories and thick subcategories that
correspond to each other under the bijection $f$; see figure \ref{fig:wideandthick}. We begin with a simple motivating example.

\begin{example} Let $R = \mathbb{Z}$, or more generally a PID. Since $R$ is clearly regular and coherent, by theorem \ref{thm:hoveystheorem} we know that 
the wide subcategories of $\Wide(R)$ and thick subcategories of $\Thick(R)$ are in bijection with the thick supports of 
$\Spec(\mathbb{Z})$.  So if $S$ is any proper thick support of $\Spec(\ints)$, then the corresponding wide subcategory is the collection of 
all finite abelian groups whose torsion is contained in  $S$. It can be easily verified that the Grothendieck group of this subcategory 
is a  free abelian group of rank $|S|$. (The $p$-length of a finite abelian group, for each $p$ in $S$, gives the universal Euler 
characteristic function.) If $S = \Spec(\mathbb{Z})$, then the corresponding wide subcategory is clearly the collection of all finitely 
generated abelian groups. The Grothendieck group of the latter is well-known to be infinite cyclic.  (The free rank of the finitely generated 
abelian group gives the universal Euler characteristic function.) It turns out that one arrives at the same answers by computing the 
Grothendieck groups for the  thick subcategories of perfect complexes over a PID; see ~\cite[Theorem 4.12]{cheb06a}.
\end{example}

So motivated by this example, we now prove that the above bijection between the lattices $\mathcal{L}_{\Wide} (R)$ 
and $\mathcal{L}_{\Thick}(R)$ preserves $K$-theory. We will need the following lemma that characterises small objects in the 
derived category of a regular coherent ring.

\begin{lemma} \cite[Note 9.8]{ch} \label{lemma:regularcoherent-small} Let $R$ be a regular coherent ring. Then a complex $X$ in $D(R)$ is small, i.e., $X$ belongs to
$\Thick(R)$, if and only if each $H_n(X)$ is finitely presented and only finitely many are non-zero.
\end{lemma}

The main theorem of this section then states:

\begin{thm}\label{thm:K_0wide} Let $R$ be a regular coherent ring. If $\C$ is any wide subcategory of $\Wide(R)$, then
\[ K_0(\C) \cong K_0(f(\C)),\]
where $f(\C)$ is the thick subcategory  $\{ X  \in \Thick(R): H_n(X) \subseteq \C \; \text{for all}\ \;n \}.$
\end{thm}

\begin{proof} Since these Grothendieck groups are quotients of free abelian groups on the appropriate isomorphism classes
of objects, we first define maps on the free abelian groups. Let $FW$ denote the the free abelian group on the isomorphism 
classes of objects in $\C$ and let $FT$ denote the  free abelian group on the isomorphism 
classes of objects in $f(\C)$. Now we define maps 
\[  FW \stk{\phi} FT \stk{\psi} FW.  \]
Since $FW$ and $FT$ are free abelian groups, it suffices to define the maps on the basis elements. So we define $\phi(M) := M[0]$ and $\psi(X) := \sum (-1)^i H_i(X)$. 
In order to show that these maps give  isomorphism in $K$-theory, we have to check the following 
things.

\begin{enumerate}
\item Both $\phi$ and $\psi$ are well-defined group homomorphisms.
\item $\phi$ descends to $K_0(\C)$ and $\psi$ descends to $K_0(f(\C))$.
\item $\phi$ and $\psi$ are inverses to each other. 
\end{enumerate}

\vskip 2mm
\noindent
(1) $M$ being a finitely presented module over a regular coherent ring, it  has finite projective dimension \cite[Theorem 6.2.1]{sarah}.
This implies \cite[Prop. 3.4]{wide} that $M[0]$ is small. Further $H_*(M[0]) = M$ belongs to $\C$, so by definition
$M[0]$ is in $f(\C)$. This shows that $\phi$ is well defined. As for $\psi$, if $X \in f(\C)$, then by definition 
$H_n(X) \in \C$ for all $n$. So $\sum (-1)^i H_i(X)$ is a well defined element of $FW$.

\vskip 2mm
\noindent
(2) To show that $\phi$ descends, all we need to show is that whenever $0 \rightarrow A \rightarrow B \rightarrow C \rightarrow 0$
is a short exact sequence in $\C$, then $\phi(B)=\phi(A) + \phi(C)$. First observe that when we replace each of these three modules with their
projective resolutions, the resolutions will be perfect complexes because the ring is coherent and regular (see the proof of \cite[Prop. 3.4]{wide}). 
So let $P_A$ and $P_C$ be projective resolutions of the modules $A$ and $C$ respectively. Then by the Horseshoe lemma \cite[Lemma 2.2.8]{wei} we know that there is a
projective resolution $P_B$ that fits into a short exact sequences (of complexes) $0\rar P_A \rar P_B \rar P_C \rar 0$. 
From this it follows that  $[P_B] = [P_A] + [P_C]$ or equivalently $[B[0]] = [A[0]] + [C[0]]$, which translates to $f(B) = f(A) + f(C)$. 
In other words, we have shown that the functor $M \mapsto M[0]$ from $\Wide(R)$ to $\Thick(R)$
sends exact sequences to exact triangles. Therefore $\phi$ descends to $K_0 (\C)$.

Now to see that the map $\psi$ descends to $K_0(f(\C))$, observe that the functor $H_*(-)$ takes exact triangles to long exact sequences.
It now follows that if  $X \rar Y \rar Z \rar \Sigma\, X$ is a triangle in $f(\C)$, then the relation 
\[ \sum (-1)^i H_i(Y) = \sum (-1)^i H_i(X) + \sum (-1)^i H_i(Z)\]
holds in $K_0(\C)$.

\vskip 2mm
\noindent
(3) Clearly $\psi \phi (M) = \psi(M[0]) = M.$ To show that the other composition is also identity, first note that
$\phi \psi (X) =\phi (\sum(-1)^i H_i(X)) = \sum (-1)^i H_i(X)[0]$. It remains to show that this
last expression is equal to $X$  in $K_0(f(\C))$. To this end,  let $l$ denote the largest integer such that $H_i(X)=0$ for all $i > l $ 
(recall that our complexes are homologically graded). Then one gets the following exact triangles in $D(R)$; see, for instance, \cite[chapter IV, section 4.5]{GelMan}, 

\begin{center}
$X_{\le l-1} \rar X \rar X_{= l} \rar \Sigma\, X_{\le l-1}$,\\
$X_{\le l-2} \rar X_{\le l-1} \rar X_{= l-1} \rar \Sigma \, X_{\le l-2} $,\\
$X_{\le l-3} \rar X_{\le l-2} \rar X_{= l-2} \rar \Sigma \, X_{\le l-3}$  \hspace{3 mm} \text{etc.}
\end{center}
(Here $X_{\le i}$ is a complex whose homology agrees with that of $X$ in dimensions $\le i$ and has no homology in dimensions
$> i$, and similarly $X_{=i}$ is a complex whose homology agrees with that of $X$ in dimension $i$ and has no homology elsewhere.)
\noindent
Now we claim that these triangles belong to  $f(\C)$. First of all, the homology 
groups of these complexes belong to $\C$, and further it is clear from the long exact sequences in homology that the homology groups of
complexes in the above exact triangles are concentrated only in a finite range. Now we invoke lemma \ref{lemma:regularcoherent-small}
to conclude that these complex are indeed small, i.e, they belong to $\Thick(R)$.

The above exact triangles give the following equations in $K_0 (f(\C))$:
\begin{eqnarray*}
[X] &=&  [X_{\le l-1}] + [X_{= l}]     \\
    &=&  [X_{\le l-2}] + [X_{= l-1}] + [X_{= l}] \\
    &=&  [X_{\le l-3}] + [X_{= l-2}] + [X_{= l-1}] + [X_{= l}] \\
    &=&  \sum_{i \le l} [X_{= i}]
\end{eqnarray*}
It is easy to see that $X_{=i}$ is quasi-isomorphic to $H_i(X)[i]$.
By replacing this quasi-isomorphism in the last equation we get 
$[X] = \sum [H_i(X)[i]] =  \sum (-1)^i [H_i(X)[0]].$
That proves that $\phi \psi (X) = X$, and we are done.
\end{proof}

\begin{cor} If $R$ is a regular coherent ring, then 
\[ K_0 ( D^b (\Wide(R) )) \cong K_0 (D^b( \proj \,R)).\]
\end{cor}

\begin{proof} It is well-known that if $\A$ is any essentially small abelian category, then
$K_0(\A) \cong K_0(D^b(\A))$; see \cite[Chapter VIII, Page 357]{sga5}. In particular, if $\A = \Wide(R)$, the abelian category of finitely
presented $R$-modules, we get 
\[K_0(\Wide(R)) \cong K_0(D^b(\Wide(R))).\]  
By theorem \ref{thm:K_0wide} we know that $K_0(\Wide(R)) \cong K_0(D^b(\proj \, R))$. So the corollary follows by combining these
two isomorphisms.
\end{proof}

\section{Decompositions for wide subcategories} \label{sec:KSwide}

In this section we study Krull-Schmidt type decompositions for wide subcategories. We denote the coproduct in our categories by the symbol $\amalg$, and
$\mathbf{0}$ will denote the trivial category in which each object is isomorphic to the zero object.

\begin{defn} A  \emph{decomposition} of a wide subcategory $\W$  of an abelian category $\A$ is 
a collection of wide subcategories $(\W_i)_{i \in I}$ in $\A$ such that, 
\begin{enumerate}
\item $\W_i \bigcap \W_j = \mathbf{0}$ for all $i \ne j$.
\item An object $M$ is in $\W$ if and only if  $M \cong  \amalg \, M_i$, where $M_i \in \W_i$ and all but finitely many of the 
$M_i$ are $0$. 
\end{enumerate}
Such a decomposition is said to be a \emph{Krull-Schmidt} decomposition, if in addition the following two conditions hold.
\begin{enumerate}
\item Each $\W_i$ is indecomposable: i.e.,
$\W_i \ne \mathbf{0}$, and if $\mathcal{P}$ and $\mathcal{Q}$ form a decomposition of $\W_i$, then 
either $\mathcal{P} = \mathbf{0} \ \text{or\;}  \mathcal{Q} = \mathbf{0}.$
\item Uniqueness: If $(\mathcal{V}_j)_{j \in J}$ is any other collection of wide subcategories that satisfy the above three properties, 
then $\W_i \cong \mathcal{V}_i$ up to a permutation of indices.
\end{enumerate}
A (Krull-Schmidt) decomposition $(\W_i)_{i \in I}$ of $\W$, if it exists, will be denoted by 
\[\W = \coprod_{i \in I} \, \W_i.\] 
Analogously, by  a \emph{thick decomposition} of a thick support $A$ of $\Spec(R)$, we will mean a decomposition $A = \bigcup A_i$ into thick supports, 
where $A_i \cap A_j = \emptyset$ if $i \ne j$. A thick decomposition of $A$ is \emph{Krull-Schmidt} if the $A_i$ are non empty, do not admit non-trivial 
thick decompositions, and if any such decomposition of $A$ is unique up to a permutation of the subsets $S_i$. 
\end{defn}

Naturally we show that the wide subcategories admit a Krull-Schmidt decomposition by showing that their corresponding thick supports admit one. We illustrate
this technique in the following prototypical example and a well-known lemma from commutative algebra

\begin{example} Let $\W$ denote the category of finite abelian groups. $\W$ is then clearly a wide subcategory of abelian groups.
The thick support corresponding to $\W$ is the set $A$ of non-zero primes in $\ints$. So the Krull-Schmidt decomposition of $A$ is given by
\[ A = \bigcup_{p \; \text{prime}} \{ (p) \}.\]
Now if we declare  $\W_p$ to be  the category of finite $p$-torsion abelian groups, it is  easily seen that each $\W_p$ is a wide subcategory of abelian groups
and that
\[ \W = \coprod_p \, \W_p\]
is the Krull-Schmidt decomposition of $\W$. Note that the decomposition of $\W$ here corresponds to that of $A$ under the bijection $\xi$ of theorem 
\ref{thm:hoveystheorem}.  
\end{example}

The next lemma, which is a well-known fact from commutative algebra, is another model on which our Krull-Schmidt decomposition is based. 

\begin{lemma} \label{cor:splittingformodules} Let $R$ be a noetherian ring and $M$ a finitely generated $R$-module. Then $M$
admits a unique splitting
\[M \cong \bigoplus_{i = 1}^n \; M_i, \]
for some $n$ such that the supports of the $M_i$ are pairwise disjoint and indecomposable.
\end{lemma}

Motivated by the above example and lemma, we now prepare towards our Krull-Schmidt decomposition theorem for wide subcategories. 

\begin{lemma} \label{lemma:ext^i=0} Let $M$ and $N$ be finitely presented modules with disjoint supports over a commutative ring $R$. Then
\[ \Ext^i_R(M , N) = 0 \ \ \text{for all} \;  \; i \;  \ge 0 .\]
\end{lemma}

\begin{proof} Recall that localisation commutes with the Ext functor for finitely presented $R$-modules.
Now since $M$ and $N$ have disjoint supports, we get 
\[ \Ext^i_R {(M, N)}_p = \Ext^i_{R_p} ({M}_p, {N}_p)  = 0 \ \ \text{for all} \; \; p \; \text{in} \; \Spec(R). \]
Therefore $\Ext^i_R(M, N) = 0$ for all $i \ge 0$.
\end{proof}

\begin{prop} \label{prop:0} Let $R$ be a regular coherent ring.  If $\W_1$ and $\W_2$ are wide subcategories of $\Wide(R)$ such that 
$\W_1 \bigcap \W_2 =\mathbf{0}$, then so is their coproduct $\W_1 \coprod \W_2$.
\end{prop}

\begin{proof} Let $A_1$ and $A_2$ denote the thick supports corresponding to the wide subcategories $\W_1$ and $\W_2$ respectively.
Since Hovey's bijection ($\xi$, see figure \ref{fig:wideandthick}) is a map of posets, it is clear that $A_1 \bigcap A_2$ is empty. So if $M_1$ and $M_2$ are modules in 
$\W_1$ and $\W_2$ respectively, then $M_1$ and $M_2$ will have disjoint supports. So by lemma \ref{lemma:ext^i=0}, every map  in the category 
$\W_1 \coprod \W_2$ splits a sum of two maps, one from $\W_1$ and the other from $\W_2$. 
This shows that $\W_1 \coprod \W_2$ is an abelian subcategory of $\Wide(R)$.

Now we argue in a similar fashion that $\W_1 \coprod \W_2$ is also closed under extensions.
Towards this, consider a short exact sequence  
\[ 0 \rar M_1 \oplus M_2 \rar M \rar M_1'\oplus M_2'\rar 0 ,\]
in $\Wide(R)$, where $M_1\oplus M_2$ and $M_1' \oplus M_2'$ belong to $\W_1 \coprod \W_2$. This short exact sequence corresponds to a class in 
$\Ext^1_R (M_1' \oplus M_2', M_1\oplus M_2)$. Since $\Ext$ preserves finite coproducts in both components, we get 
\[\Ext^1_R (M_1' \oplus M_2', M_1\oplus M_2) \cong \underset{i,j}{\bigoplus} \Ext^1_R(M_i', M_j )  .\]
Since the supports of $M_i'$ and $M_j$ are disjoint, we know from lemma \ref{lemma:ext^i=0} that  $\Ext^1_R(M_i', M_j)$ is zero for $i \ne j$.
Therefore the last equation of $\Ext$ groups now simplifies to 
\[\Ext^1_R (M_1' \oplus M_2', M_1\oplus M_2) = \Ext^1_R(M_1', M_1) \oplus \Ext^1_R(M_2' ,M_2) .  \]
So the upshot is that every extension of $M_1 \oplus M_2$ and $M_1' \oplus M_2'$ is a sum of two extensions: one obtained from an extension 
of $M_1$ and $M_1'$ and the other from $M_2$ and $M_2'$. This shows that  $\W_1 \coprod \W_2$ is also closed under extensions. 
\end{proof}

Now we establish the strong connection between decompositions of thick supports and those of wide subcategories.
If $A$ is a thick support of $\Spec(R)$ we denote the corresponding wide subcategory  by $\W_A$.

\begin{prop} \label{prop:1} Let $R$ be a regular coherent ring. If $A = A_1 \bigcup A_2$ is a thick decomposition of a thick support $A$ in $\Spec(R)$, 
then this induces a decomposition of the associated wide subcategories in $\Wide(R)$:
 \[ \W_A = \W_{A_1} \coprod \W_{A_2}. \]
\end{prop}

\begin{proof} Clearly $\W_{A_1} \bigcap \W_{A_2} = \mathbf{0}$ because $A_1$ and $A_2$ are disjoint. We only have to show that the objects in $\W_A$ are the 
coproducts of objects in $\W_{A_1}$ and $\W_{A_2}$. This is shown in an indirect way as follows.
By proposition \ref{prop:0} we know that  $\W_{A_1} \coprod \W_{A_2}$ is a wide subcategory. So we will be done if we can show that the thick 
support corresponding to   $\W_{A_1} \coprod \W_{A_2}$ is $A$. This follows when we observe that the union of the
supports of the  modules in  $\W_{A_1} \coprod \W_{A_2}$ is clearly $A_1 \bigcup A_2 (= A)$. Hence  $\W_{A_1} \coprod \W_{A_2} = \W_A$.
\end{proof}

\begin{prop} \label{prop:2} Let $R$ be a regular coherent ring and let $\W$ be a wide subcategory of $\Wide(R)$ corresponding to a thick support $A$. 
If  $\W = \W_1 \coprod \W_2$ is a decomposition of $\W$ into wide subcategories, then this induces a decomposition $A = A_1 \bigcup A_2$ of the thick support $A$,
where $A_i$ is the thick support corresponding to $\W_i$. 
\end{prop}

\begin{proof} Again using the fact that Hovey's bijection is a map of posets, it is clear that $A_1$ and $A_2$ are disjoint thick supports of $\Spec(R)$ that are 
contained in $A$. To see that we have equality ($A = A_1 \bigcup A_2$) observe that each module $M$ in $\W$ splits as $M = M_1 \oplus M_2$, and recall
that $\Supp(M) = \Supp(M_1) \bigcup \Supp(M_2).$ 
\end{proof}

Combining propositions \ref{prop:1} and \ref{prop:2}, we arrive at the following decomposition result.

\begin{thm} \label{thm:KSforregularcoherent} Let $R$ be a regular coherent ring. A wide subcategory of $\Wide(R)$ admits Krull-Schmidt decomposition if and 
only if the  corresponding thick support admits one. 
\end{thm}

Now the question that has to be addressed is whether every thick support of $\Spec(R)$ admits a Krull-Schmidt decomposition. We show that 
this is possible if we assume the ring to be noetherian.

\begin{prop} \label{prop:graph} Let $R$ be a noetherian ring and let $A$ be a thick support of $\Spec(R)$. 
Then there exists a Krull-Schmidt decomposition $\bigcup A_i$ for $A$.
\end{prop}

\begin{proof} 
It is well-known that the set of prime ideals in a noetherian ring satisfies the descending
chain condition \cite[Corollary 11.12]{am}.  To start, let $A$ be a thick support in $\Spec(R)$ and let $(p_i)_{i \in I}$ be the collection 
of all minimal elements in $A$ -- i.e., primes $p$ in $A$ which do not contain any other prime in $A$.  It is now clear (using the above fact about noetherian rings)
that every prime $p \in A$ contains a minimal element $p_i$ in $A$, therefore $A = \bigcup_{i \in I} V(p_i)$. (Also note that each 
$V(p_i)$ is a closed subset of $\Spec(R)$ and hence a thick support.)  Now define a graph $G_A$ of $A$ as follows:
The  vertices are the minimal primes $(p_i)_{i \in I}$ in $A$, and two vertices $p_i$ and $p_j$ are adjacent if and only if 
$V(p_i) \cap V(p_j) \ne \emptyset$. 
Let $(C_k)_{k \in K}$ be the connected components of this graph and for each $C_k$ define a thick support
\[A_k:= \underset{p_i \in C_k}{\bigcup}{ V(p_i)}.\]
By construction it is clear that $\bigcup A_k$ is a thick decomposition of $A$. It is not hard to see that each $A_k$ is indecomposable.
Finally the uniqueness part: Suppose $\bigcup B_k$ is a Krull-Schmidt decomposition of $A$. 
It can be easily verified that the minimal primes in $B_k$ are precisely the minimal primes of $A$ that are contained in $B_k$. Thus the 
Krull-Schmidt decomposition $\bigcup B_k$ gives a partition of the set of minimal primes in $A$. This partition induces a decomposition of 
$G_A$ into its connected (since each $B_k$ is indecomposable) components. Since the decomposition of a graph into its connected 
components is unique, the uniqueness of Krull-Schmidt decomposition follows.
\end{proof}

The main theorem of this section then states:

\begin{thm} \label{thm:KSwide} Let $R$ be a noetherian regular ring. Then every wide subcategory of $\Wide(R)$ admits a Krull-Schmidt
decomposition. Conversely, given indecomposable wide subcategories  $\W_i$  of $\Wide(R)$ such that $\W_i \bigcap \W_j = \mathbf{0}$ for $i \ne j$, there 
exists a unique wide subcategory $\W \subseteq \Wide(R)$ such that $\W = \coprod_{i \in I} \W_i$ is a Krull-Schmidt decomposition for $\W$.
\end{thm}

\begin{proof} The first part of the theorem follows from theorem \ref{thm:KSforregularcoherent} and proposition \ref{prop:graph}. For the converse,
define $\W$ to the full subcategory of all finite coproducts of objects from the $\W_i$. It follows from proposition \ref{prop:0} that $\W$ thus defined is wide.
\end{proof}

We now give some corollaries of our decomposition theorem. We begin by observing that K-theory is preserved under these
decompositions.

\begin{cor} Let $R$ be a noetherian regular ring and let $\W$ be a wide subcategory of finitely generated $R$ modules. Then
\[K_0(\W) \cong \underset{i \in I}{\bigoplus} \; K_0(\W_i),\]
where $\coprod_{i \in I} \W_i$ is the Krull-Schmidt decomposition of $\W$.
\end{cor}

\begin{proof} Consider the Krull-Schmidt decomposition $\coprod_{i \in I} \W_i$ of $\W$. We know from lemma \ref{lemma:ext^i=0} that for 
$i \ne j$, $\Hom(\W_i, \W_j) = 0$. So it is clear that every short exact sequence in $\W$ breaks (uniquely up to isomorphism and permutation) 
as a sum of finitely many short exact sequences in the categories $\W_i$. It  follows that $K_0(\W_i) \cong \oplus_{i \in I} K_0(\W_i)$.
\end{proof}

The next corollary is a categorical  characterisation of local rings amongst noetherian regular rings.

\begin{cor} A noetherian regular ring $R$ is local if and only if every wide subcategory of $\Wide(R)$ is indecomposable.
\end{cor}

\begin{proof} If $R$ is local then it is clear that every thick support of $\Spec(R)$ contains the unique maximal ideal. In particular,
we cannot have two non-empty disjoint thick supports and therefore every wide subcategory of $\Wide(R)$ must be indecomposable. Conversely, if $R$
is not local, then the wide subcategory supported on the set of maximal ideals is decomposable. So we are done.
\end{proof}

\begin{rem} It is worth noting that our proofs in this section give the following stronger statements.
\begin{enumerate}
\item Let $R$ be a coherent ring  for which the map (see theorem \ref{thm:hoveystheorem})
\[f: \mathcal{L}_{\Wide}(R) \rar \mathcal{L}_{\Thick}(R)\]
is an isomorphism. Then a wide subcategory of $\Wide(R)$ admits Krull-Schmidt decomposition if and 
only if the  corresponding thick support admits one. 

\item Let $R$ be a regular coherent ring with the following properties.

(a) Every open subset of $\Spec(R)$ is compact.

(b) $\Spec(R)$ satisfies the descending chain condition.

\noindent
Then every wide subcategory of finitely presented modules admits a Krull-Schmidt decomposition. 
\end{enumerate}

However plausible these decompositions might sound, it is not clear how one would arrive at them
without using Hovey's classification. It would be interesting to find a  direct way to get these decompositions. Such an approach might 
actually give us a much stronger result, possibly without the regularity assumption on the ring.

For the curious reader we mention that these decompositions have also been studied in the context of triangulated categories. 
\cite{Kr} and \cite{KS} discuss Krull-Schmidt decompositions for thick subcategories 
of small objects in triangulated categories such as stable module categories over group algebras, derived categories of rings, and the stable homotopy category of spectra.
\end{rem}

\bibliographystyle{alpha}

\end{document}